\documentclass[10pt,a4paper]{article}
\begin{document}
\begin{center}
 {\huge\textbf{Stationary Solutions of the Klein-Gordon Equation in a Potential Field}}
 \vspace{0.5cm}

{\large\texttt{Guangqing Bi}}

{\normalsize\texttt{College of Physics and Electronic Information,
Yan'an University, Yan'an, Shaanxi, China, 716000}}

{\normalsize\texttt{yananbiguangqing@sohu.com}} \vspace{0.3cm}

{\large\texttt{Yuekai Bi}}

{\normalsize\texttt{School of Electronic and Information
Engineering, BUAA, Beijing, China, 100191}}

{\normalsize\texttt{yuekaifly@163.com}}
\end{center}
 \noindent

\begin{abstract}
We seek to introduce a mathematical method to derive the
Klein-Gordon equation and a set of relevant laws strictly, which
combines the relativistic wave functions in two inertial frames of
reference. If we define the stationary state wave functions as
special solutions like
$\Psi(\mathbf{r},t)=\psi(\mathbf{r})e^{-iEt/\hbar}$, and define
$m=E/c^2$, which is called the mass of the system, then the
Klein-Gordon equation can clearly be expressed in a better form when
compared with the non-relativistic limit, which not only allows us
to transplant the solving approach of the Schr\"{o}dinger equation
into the relativistic wave equations, but also proves that the
stationary solutions of the Klein-Gordon equation in a potential
field have the probability significance. For comparison, we have
also discussed the Dirac equation. By introducing the concept of
system mass into the Klein-Gordon equation with the scalar and
vector potentials, we prove that if the Schr\"{o}dinger equation in
a certain potential field can be solved exactly, then under the
condition that the scalar and vector potentials are equal, the
Klein-Gordon equation in the same potential field can also be solved
exactly by using the same method.
\end{abstract}

\textbf{Keywords:} Relativistic wave equations; Klein-Gordon
equation; Bound state

\textbf{PACS:} 03.65.Pm, 03.65.Ge

\section{Introduction}

Both special relativity and experiments indicate that, the mass of a
many-particle system in a bound state is less than the sum of the
rest mass of every particle forming the system, and the difference
gives the mass defect of the system, while the product of the mass
defect and the square of the speed of light gives the binding energy
of the system. As the binding energy is quantized, the sum of it and
the rest mass of every particle forming the system is the energy
level of the system. For instance, the mass of an atomic nucleus is
obviously less than the sum of the rest mass of every nucleon
forming the atomic nucleus. Therefore, in order to express the mass
defect explicitly, there is a necessity to introduce the concept of
system mass, which differs from the sum of the rest mass of every
particle forming the system. By introducing this concept, we can
express relativistic wave equations in a better form when compared
with the non-relativistic limit. By means of this method, we are
able to solve the Klein-Gordon equation just like solving the
Schr\"{o}dinger equation, and the solutions can also have the
probability significance just like the solutions of the
Schr\"{o}dinger equation. As the Schr\"{o}dinger equation has been
thoroughly studied, it surely creates conditions for studying the
Klein-Gordon equation thoroughly.

\section{Relativistic Wave Equations and the Probability Significance of Their Stationary Solutions}

As we know, for a free particle, its energy and momentum are both
constant, therefore, it is quite natural to assume that its matter
wave is a plane wave. According to the de Broglie relation
\[\hbar\mathbf{k}=\mathbf{p},\qquad{E}=\hbar\omega.\]
We have the wave function of a free particle:
\begin{equation}\label{1}
\Psi_p(\mathbf{r},t)=Ae^{-\frac{i}{\hbar}(Et-\mathbf{p\cdot{r}})}.
\end{equation}
Where $E$ is the total energy of the particle containing the
intrinsic energy $m_0c^2$, and $\mathbf{p}$ is the momentum of the
particle.

The quantization method of quantum mechanics is assuming $E$ and
$\mathbf{p}$ are correspondingly equivalent to the following two
differential operators:
\[E\rightarrow i\hbar\frac{\partial}{\partial{t}},\qquad\mathbf{p}\rightarrow{-i}\hbar\nabla.\]
This relation is obviously tenable for the wave function of a free
particle, while arbitrary wave function is equal to the linear
superposition of the plane waves of free particles with all possible
momentum, i.e.
\begin{eqnarray}\label{1a}
\Psi(\mathbf{r},t)&=&\frac{1}{(2\pi\hbar)^{3/2}}\int\!\!\!\int\limits^\infty_{-\infty}\!\!\!\int{c}(\mathbf{p},t)e^{-\frac{i}{\hbar}(Et-\mathbf{p\cdot{r}})}\,dp_xdp_ydp_z,\\
c(\mathbf{p},t)&=&\frac{1}{(2\pi\hbar)^{3/2}}\int\!\!\!\int\limits^\infty_{-\infty}\!\!\!\int\Psi(\mathbf{r},t)e^{\frac{i}{\hbar}(Et-\mathbf{p\cdot{r}})}\,dx\,dy\,dz.
\end{eqnarray}
Therefore, this quantization method is also tenable for arbitrary
wave functions. Later, we will prove the right-hand side of
(\ref{1}) is relativistic, thus this quantization method itself is
relativistic. Then why the Schr\"{o}dinger equation derived from
this method is non-relativistic? That is because the relation
between $E$ and $\mathbf{p}$ is non-relativistic, which is due to
the fact that the kinetic energy is non-relativistic. Therefore, if
we want to establish the relativistic wave equations, we need to
introduce the relativistic kinetic expression.

Assuming that any particle with the rest mass $m_0$, no matter how
high the speed is, no matter it is in a potential field or in free
space, and no matter how it interacts with other particles, its
kinetic energy is:
\begin{eqnarray}\label{2a}
E_k=(c^2p^2+m_0{}^2c^4)^{1/2}-m_0c^2.
\end{eqnarray}

If $E_k+m_0c^2$ is denoted by $E$, then (\ref{2a}) can be expressed
as
\[E^2=c^2p^2+m_0{}^2c^4.\]
Thus, it leaves us a matter of mathematical techniques. Therefore,
we introduce the mathematical method in reference
{\cite{bi97}-\cite{bi10}}. Power functions and exponential functions
play special roles in this method, which are called the base
functions as we can establish mapping relations between them and
arbitrary functions in a certain range. For instance, in quantum
mechanics, (\ref{1a}) determines the mapping relations between wave
functions of free particles and arbitrary wave functions. Similarly
to reference {\cite{bi97}-\cite{bi10}}, we can introduce the base
function and relevant concepts in quantum mechanics:

\textbf{Definition 1} The right-hand side of (\ref{1}) is defined as
the base function of quantum mechanics, where $E$ and $\mathbf{p}$
are called the \textbf{characters} of base functions, while $E$ and
$\mathbf{p}$ are not only suitable for free particles, but also
suitable for any system, and the relation between $E$ and
$\mathbf{p}$ is called the characteristic equation of wave
equations. Different system has different characteristic equations.
For instance, $E^2=c^2p^2+m_0{}^2c^4$ is the characteristic equation
for free-particle system.

According to differential laws, we have
\begin{equation}\label{2}
i\hbar\frac{\partial}{\partial{t}}\Psi_p=E\Psi_p,\quad{-i}\hbar\nabla\Psi_p=\mathbf{p}\Psi_p.
\end{equation}

\textbf{Definition 2} Let $m_0=m_{01}+m_{02}+\cdots+m_{0N}$ be the
total rest mass of an $N$-particle system, $E'$ be the sum of the
kinetic energy and potential energy of all the $N$ particles, then
the actual mass of the system, which is called the \textbf{system
mass}, is defined as
\begin{equation}\label{3}
m=m_0+\frac{1}{c^2}E'.
\end{equation}

\textbf{Definition 3}  If the system is in a bound state($E'<0$),
then the absolute value of $E'$ is
\[|E'|=m_0c^2-mc^2=\triangle{m}c^2,\] which is called the binding
energy of the system, where $\triangle{m}=m_0-m$ is the mass defect
of the system.

\textbf{Definition 4} The total energy of the system $E$ is defined
as the sum of the rest energy, kinetic energy and potential energy
of all the particles forming the system, i.e. $E=m_0c^2+E'$.

According to Definition 2 and 4, the total energy of the system is
equal to the product of the system mass and the square of the speed
of light, i.e. $E=mc^2$, thus the system mass is uniquely determined
by the energy level of the system.

\textbf{Definition 5} In relativistic quantum mechanics, the
stationary state wave function is defined as the following special
solution, i.e.
\begin{equation}\label{4'}
\Psi(\mathbf{r},t)=\psi(\mathbf{r})e^{-iEt/\hbar}.
\end{equation}

Let $U(\mathbf{r})$ be the potential energy of a particle moving in
an potential field, then according to Definition 2, we have
$E_k=E'-U$. Combining it with (\ref{2a}), we have
\begin{equation}\label{4}
    E'-U=(c^2p^2+m_0{}^2c^4)^{1/2}-m_0c^2.
\end{equation}
This is the characteristic equation of relativistic wave equations
for a single-particle system. In order to make it easier to solve
the corresponding relativistic wave equation, the characteristic
equation (\ref{4}) should be transformed to remove the fractional
power, then we have
\begin{equation}\label{5}
(E'-U+m_0c^2)^2=c^2p^2+m_0{}^2c^4.
\end{equation}

Expanding the left-hand side of (\ref{5}) and by using (\ref{3}), we
have
\begin{eqnarray}\label{6}
E'&=&\frac{p^2}{m_0+m}+\frac{2m}{m_0+m}U-\frac{U^2}{(m_0+m)c^2},\nonumber\\
E&=&\frac{p^2}{m_0+m}+\frac{2m}{m_0+m}U-\frac{U^2}{(m_0+m)c^2}+m_0c^2.
\end{eqnarray}
The essence of this expression is the relativistic Hamiltonian.
Therefore, taking (\ref{6}) as the characteristic equation,
multiplying both sides of the equation by the base function
$\Psi_p(\mathbf{r},t)$, and by using (\ref{2}), we have
\[i\hbar\frac{\partial\Psi_p}{\partial{t}}=-\frac{\hbar^2}{m_0+m}\nabla^2\Psi_p+\frac{2m}{m_0+m}U\Psi_p-\frac{U^2}{(m_0+m)c^2}\Psi_p+m_0c^2\Psi_p.\]
According to (\ref{1a}), in the operator equation which is tenable
for the base function $\Psi_p$, as long as each operator in the
operator equation is a linear operator and each linear operator does
not explicitly contain the characters $E$ and $\mathbf{p}$ of
$\Psi_p$, then this operator equation is also tenable for an
arbitrary wave function $\Psi(\mathbf{r},t)$. Whereas, considering
that the system mass $m$ is equivalent to the character $E$ of
$\Psi_p$, this operator equation is not tenable for arbitrary wave
functions, but tenable for an stationary state wave function like
(\ref{4'}). Thus we have:

A particle with the mass $m_0$ moving in the potential field
$U(\mathbf{r})$ can be described by the wave function
$\Psi(\mathbf{r},t)$, an arbitrary stationary state wave function
$\Psi$ satisfies the following relativistic wave equation
\begin{equation}\label{7}
i\hbar\frac{\partial\Psi}{\partial{t}}=-\frac{\hbar^2}{m_0+m}\nabla^2\Psi+\frac{2m}{m_0+m}U\Psi-\frac{U^2}{(m_0+m)c^2}\Psi+m_0c^2\Psi.
\end{equation}
The corresponding $\psi(\mathbf{r})$ satisfies the following
relativistic wave equation
\begin{equation}\label{8}
E\psi=-\frac{\hbar^2}{m_0+m}\nabla^2\psi+\frac{2m}{m_0+m}U\psi-\frac{U^2}{(m_0+m)c^2}\psi+m_0c^2\psi.
\end{equation}
Or, if the energy corresponding to the energy operator is
interpreted as $E'$ in Definition 2, then the relativistic wave
equation is
\begin{equation}\label{9}
i\hbar\frac{\partial\Psi}{\partial{t}}=-\frac{\hbar^2}{m_0+m}\nabla^2\Psi+\frac{2m}{m_0+m}U\Psi-\frac{U^2}{(m_0+m)c^2}\Psi.
\end{equation}
The corresponding stationary form is
\begin{equation}\label{10}
E'\psi=-\frac{\hbar^2}{m_0+m}\nabla^2\psi+\frac{2m}{m_0+m}U\psi-\frac{U^2}{(m_0+m)c^2}\psi.
\end{equation}
As $E-m_0c^2=E'$, (\ref{7}) and (\ref{9}) are equivalent. In the
non-relativistic limit, we have
\[m\rightarrow{m_0},\quad\frac{U^2}{(m_0+m)c^2}\rightarrow0,\]
(\ref{9}) approaches the Schr\"{o}dinger equation.

Let $w(\mathbf{r},t)=\Psi^*(\mathbf{r},t)\Psi(\mathbf{r},t)$,
differentiating both sides of the equation with respect to time $t$,
and considering
\[\frac{\partial\Psi}{\partial{t}}=\frac{i\hbar}{m_0+m}\nabla^2\Psi+\frac{1}{i\hbar}\frac{2m}{m_0+m}U\Psi-\frac{1}{i\hbar}\frac{U^2}{(m_0+m)c^2}\Psi+\frac{1}{i\hbar}m_0c^2\Psi,\]
\[\frac{\partial\Psi^*}{\partial{t}}=-\frac{i\hbar}{m_0+m}\nabla^2\Psi^*-\frac{1}{i\hbar}\frac{2m}{m_0+m}U\Psi^*+\frac{1}{i\hbar}\frac{U^2}{(m_0+m)c^2}\Psi^*-\frac{1}{i\hbar}m_0c^2\Psi^*.\]
We have the probability conservation equation, i.e.
\[\frac{\partial{w}}{\partial{t}}=\frac{i\hbar}{m_0+m}(\Psi^*\nabla^2\Psi-\Psi\nabla^2\Psi^*)=\frac{i\hbar}{m_0+m}\nabla\cdot(\Psi^*\nabla\Psi-\Psi\nabla\Psi^*)=-\nabla\cdot\mathbf{J}.\]
In other words,
$w(\mathbf{r},t)=\Psi^*(\mathbf{r},t)\Psi(\mathbf{r},t)$ can be
interpreted as the probability density. Thus we have

Let $w(\mathbf{r},t)=\Psi^*(\mathbf{r},t)\Psi(\mathbf{r},t)$ be the
probability density, then the corresponding probability current
density vector is
\begin{eqnarray}\label{15}
\mathbf{J}=\frac{i\hbar}{m_0+m}(\Psi\nabla\Psi^*-\Psi^*\nabla\Psi),
\end{eqnarray}
making the following equation tenable, i.e.
\begin{eqnarray}\label{16}
\frac{\partial{w}}{\partial{t}}+\nabla\cdot\mathbf{J}=0.
\end{eqnarray}
Here $\Psi(\mathbf{r},t)$ is a stationary state wave function like
(\ref{4'}), determined by wave function (\ref{7}) or (\ref{9}) and
natural boundary conditions.

Now, Let us introduce the concept of system mass into the Dirac
equation, so as to compare it with (\ref{10}). According to
reference {\cite{zhou}}, we can express the Dirac equation as the
following two equations:
\[c\widehat{\sigma}\cdot\mathbf{p}\psi_2+\mu_0c^2\psi_1+U\psi_1=E\psi_1,\]
\[c\widehat{\sigma}\cdot\mathbf{p}\psi_1-\mu_0c^2\psi_2+U\psi_2=E\psi_2.\]
Considering $E=\mu c^2$, solve the second equation for $\psi_2$,
then substitute the result into the first equation, thus we have
\[(\widehat{\sigma}\cdot\mathbf{p})\frac{c^2}{(\mu_0+\mu)c^2-U}(\widehat{\sigma}\cdot\mathbf{p})\psi_1+\mu_0c^2\psi_1+U\psi_1=E\psi_1.\]
Where
\begin{eqnarray*}
   & & (\widehat{\sigma}\cdot\mathbf{p})\frac{c^2}{(\mu_0+\mu)c^2-U}(\widehat{\sigma}\cdot\mathbf{p})\psi_1 \\
   &=& \frac{c^2}{(\mu_0+\mu)c^2-U}(\widehat{\sigma}\cdot\mathbf{p})(\widehat{\sigma}\cdot\mathbf{p})\psi_1
   +\left((\widehat{\sigma}\cdot\mathbf{p})\frac{c^2}{(\mu_0+\mu)c^2-U}\right)(\widehat{\sigma}\cdot\mathbf{p})\psi_1.
\end{eqnarray*}
As
$(\widehat{\sigma}\cdot\mathbf{p})(\widehat{\sigma}\cdot\mathbf{p})=p^2$,
$E=E'+\mu_0c^2$, the equation can be transformed into the following
form:
\[E'\psi_1=\frac{c^2p^2}{(\mu_0+\mu)c^2-U}\psi_1+U\psi_1+
\left((\widehat{\sigma}\cdot\mathbf{p})\frac{c^2}{(\mu_0+\mu)c^2-U}\right)(\widehat{\sigma}\cdot\mathbf{p})\psi_1.\]
As it is known, for any scalar function $U$, we have
\[(\widehat{\sigma}\cdot\mathbf{p}U)(\widehat{\sigma}\cdot\mathbf{p})\psi_1=\left\{-i\hbar(\nabla U\cdot\mathbf{p})+\hbar\widehat{\sigma}[(\nabla
U)\times\mathbf{p}]\right\}\psi_1.\] If $U$ is simply a function of
$r$, then
\[\nabla U=\frac{1}{r}\frac{dU}{dr}\mathbf{r},\]
\[-i\hbar(\nabla{U}\cdot\mathbf{p})\psi_1=-\hbar^2\nabla{U}\cdot\nabla\psi_1=-\hbar^2\frac{dU}{dr}\frac{\partial\psi_1}{\partial{r}},\]
\[(\nabla{U})\times\mathbf{p}=\frac{1}{r}\frac{dU}{dr}(\mathbf{r}\times\mathbf{p})=\frac{1}{r}\frac{dU}{dr}\mathbf{L}.\]
Thus the equation can be expressed as
\begin{eqnarray*}
  E'\psi_1 &=& \frac{c^2p^2}{(\mu_0+\mu)c^2-U}\psi_1+U\psi_1\\
   & & +\,\frac{c^2}{[(\mu_0+\mu)c^2-U]^2}\frac{dU}{dr}\left(\frac{2}{r}\mathbf{S}\cdot\mathbf{L}\psi_1-\hbar^2\frac{\partial\psi_1}{\partial{r}}\right).
\end{eqnarray*}
Where $\mathbf{S}=(1/2)\hbar\widehat{\sigma}$ is the spin angular
momentum operator. Multiplying both sides of the equation by
$(\mu_0+\mu)c^2-U$, noting that $E'=E-\mu_0c^2=\mu{c^2}-\mu_0c^2$,
we have
\begin{eqnarray*}
  E'(\mu_0+\mu)c^2\psi_1 &=& c^2p^2\psi_1+2\mu{c^2}U\psi_1-U^2\psi_1\\
   & & +\,\frac{c^2}{(\mu_0+\mu)c^2-U}\frac{dU}{dr}\left(\frac{2}{r}\mathbf{S}\cdot\mathbf{L}\psi_1-\hbar^2\frac{\partial\psi_1}{\partial{r}}\right).
\end{eqnarray*}

Therefore, the Dirac equation can be expressed as the following
form:

\parbox{11cm}{\begin{eqnarray*}\label{21}
E'\psi&=&-\,\frac{\hbar^2}{\mu_0+\mu}\nabla^2\psi+\frac{2\mu}{\mu_0+\mu}U\psi-\frac{U^2}{(\mu_0+\mu)c^2}\psi\\
&&+\,\frac{1}{(\mu_0+\mu)[(\mu_0+\mu)c^2-U]}\frac{dU}{dr}\left(\frac{2}{r}\mathbf{S}\cdot\mathbf{L}\psi-\hbar^2\frac{\partial\psi}{\partial{r}}\right).
\end{eqnarray*}}\hfill\parbox{1cm}{\begin{eqnarray}\end{eqnarray}}
Where $\mu$ is the system mass corresponding to $\mu_0$,
$\mathbf{S}$ is the spin angular momentum operator, $\mathbf{L}$ is
the orbital angular momentum operator, and $U$ is simply a function
of $r$.

\section{The System Mass of Pionic Hydrogen Atoms and Relativistic Wave Functions}

A pionic hydrogen atom is a system formed by a negative pion and an
atomic nucleus. Reference \cite{w} expounded the relativistic energy
levels of such a system. Let us take the pionic hydrogen atom for
example. We can transplant the mathematical methods of the
non-relativistic quantum mechanics in reference \cite{zhou1} into
the relativistic quantum mechanics, thus determine the system mass
and relativistic wave functions of such a system.

Considering a negative pion $\pi^-$ with the rest mass $m_0$ moving
in an nuclear electric field, taking the atomic nucleus as the
origin of coordinates, then the potential energy of the system is
$U=-Ze_s^2/r,\; e_s=e(4\pi\varepsilon_0)^{-1/2}$, according to
(\ref{10}), the relativistic stationary state wave equation of the
system is
\[E'\psi=-\frac{\hbar^2}{m_0+m}\nabla^2\psi-\frac{2m}{m_0+m}\frac{Ze_s{}^2}{r}\psi-\frac{1}{(m_0+m)c^2}\frac{Z^2e_s{}^4}{r^2}\psi.\]

The equation expressed in spherical polar coordinates is
\begin{eqnarray*}
E'\psi &=&
-\frac{\hbar^2}{m_0+m}\left[\frac{\partial}{\partial{r}}\left(r^2\frac{\partial}{\partial{r}}\right)+
\frac{1}{\sin\theta}\frac{\partial}{\partial\theta}\left(\sin\theta\frac{\partial}{\partial\theta}\right)+
\frac{1}{\sin^2\theta}\frac{\partial^2}{\partial\varphi^2}\right]\psi \\
       & & -\frac{2m}{m_0+m}\frac{Ze_s{}^2}{r}\psi-\frac{1}{(m_0+m)c^2}\frac{Z^2e_s{}^4}{r^2}\psi.
\end{eqnarray*}

By following the solving procedure of the Schr\"{o}dinger equation
in reference \cite{zhou1}, we use the method of separation of
variables to solve this equation. Let
$\psi(r,\theta,\varphi)=R(r)Y(\theta,\varphi)$, by substituting it
into the equation, we have
\begin{eqnarray*}
  \frac{(m_0+m)E'r^2}{\hbar^2} &+& \frac{1}{R}\frac{\partial}{\partial{r}}\left(r^2\frac{\partial{R}}{\partial{r}}\right)
  +2m\frac{Ze_s{}^2}{\hbar^2}r+\frac{Z^2e_s{}^4}{\hbar^2c^2} \\
    &=& -\frac{1}{Y}\left[\frac{1}{\sin\theta}\frac{\partial}{\partial\theta}\left(\sin\theta\frac{\partial{Y}}{\partial\theta}\right)
    +\frac{1}{\sin^2\theta}\frac{\partial^2Y}{\partial\varphi^2}\right]=\lambda.
\end{eqnarray*}
\begin{equation}\label{a3}
\frac{1}{r^2}\frac{\partial}{\partial{r}}\left(r^2\frac{\partial{R}}{\partial{r}}\right)+\frac{(m_0+m)E'}{\hbar^2}R
+\left[\frac{2m}{\hbar^2}\frac{Ze_s{}^2}{r}+\left(\frac{Z^2e_s{}^4}{\hbar^2c^2}-\lambda\right)\frac{1}{r^2}\right]R=0.
\end{equation}
\begin{equation}\label{a4}
\frac{1}{\sin\theta}\frac{\partial}{\partial\theta}\left(\sin\theta\frac{\partial{Y}}{\partial\theta}\right)
+\frac{1}{\sin^2\theta}\frac{\partial^2Y}{\partial\varphi^2}+\lambda{Y}=0.
\end{equation}
According to (\ref{a4}), denoting $\lambda=l(l+1),\;l=0,1,2,\ldots$,
obviously the solution of the equation is the spherical harmonics
$Y_{lm}(\theta,\varphi)$.

Now let us solve the radial equation (\ref{a3}), discussing the
situation of the bound state ($E'<0$). Let
\begin{equation}\label{a5}
    \alpha'=\left[\frac{4(m_0+m)|E'|}{\hbar^2}\right]^{1/2},\quad\beta=\frac{2mZe_s{}^2}{\alpha'\hbar^2}.
\end{equation}

Let $R(r)=u(r)/r$, considering
\[\frac{1}{r^2}\frac{d}{dr}\left(r^2\frac{dR}{dr}\right)=\frac{1}{r}\frac{d^2}{dr^2}(rR),\]
and by using the variable substitution $\rho=\alpha'r$, then the
equation (\ref{a3}) can be expressed as ($\alpha$ is the fine
structure constant, which is absolutely different from $\alpha'$):
\begin{equation}\label{a6}
    \frac{d^2u}{d\rho^2}+\left[\frac{\beta}{\rho}-\frac{1}{4}-\frac{l(l+1)-Z^2\alpha^2}{\rho^2}\right]u=0.
\end{equation}

Firstly, let us study the asymptotic behaviour of this equation,
when $\rho\rightarrow\infty$, the equation can be transformed into
the following form:
\[\frac{d^2u}{d\rho^2}-\frac{1}{4}u=0,\quad{u}(\rho)=e^{\pm\rho/2}.\]
As $e^{\rho/2}$ is in conflict with the finite conditions of wave
functions, we substitute $u(\rho)=e^{-\rho/2}f(\rho)$ into the
equation, then we have the equation satisfied by $f(\rho)$:
\begin{equation}\label{a7}
\frac{d^2f}{d\rho^2}-\frac{df}{d\rho}+\left[\frac{\beta}{\rho}-\frac{l(l+1)-Z^2\alpha^2}{\rho^2}\right]f=0.
\end{equation}
Now let us solve this equation for series solutions. Let
\begin{equation}\label{a8}
    f(\rho)=\sum^\infty_{\nu=0}b_\nu\rho^{s+\nu},\quad{b_0}\neq0.
\end{equation}
In order to guarantee the finiteness of $R=u/r$ at $r=0$, $s$ should
be no less than $1$. By substituting (\ref{a8}) into (\ref{a7}), as
the coefficient of $\rho^{s+\nu-1}$ is equal to zero, we have the
relation satisfied by $b_\nu$:
\begin{equation}\label{a9}
    b_{\nu+1}=\frac{s+\nu-\beta}{(s+\nu)(s+\nu+1)-l(l+1)+Z^2\alpha^2}b_\nu.
\end{equation}
If the series are infinite series, then when $\nu\rightarrow\infty$
we have $b_{\nu+1}/b_\nu\rightarrow1/\nu$. Therefore, when
$\rho\rightarrow\infty$, the behaviour of the series is the same as
that of $e^\rho$, then we have
\[R=\frac{\alpha'}{\rho}u(\rho)=\frac{\alpha'}{\rho}\,e^{-\rho/2}f(\rho).\]
Where $f(\rho)$ tends to infinity when $\rho\rightarrow\infty$,
which is in conflict with the finite conditions of wave functions.
Therefore, the series should only have finite terms. Let
$b_{n_r}\rho^{s+n_r}$ be the highest-order term, then $b_{n_r+1}=0$,
by substituting $\nu=n_r$ into (\ref{a9}) we have $\beta=n_r+s$. On
the other hand, the series starts from $\nu=0$, and doesn't have the
term $\nu=-1$, therefore, $b_{-1}=0$. Substituting $\nu=-1$ into
(\ref{a9}), considering $b_0\neq0$, we have
$s(s-1)=l(l+1)-Z^2\alpha^2$. Denoting $n=n_r+l+1$, then the
following set of equations can be solved for $s$ and $\beta$:
\begin{equation}\label{a10}
    \left\{\begin{array}{l}
    s(s-1)=l(l+1)-Z^2\alpha^2\\
    \beta=n_r+s\\
    n=n_r+l+1.
    \end{array}\right.
\end{equation}
We have $s=1/2\pm\sqrt{(l+1/2)^2-Z^2\alpha^2}$, taking
$s=1/2+\sqrt{(l+1/2)^2-Z^2\alpha^2}$, then
\[\beta=n_r+s=n-l-1/2+\sqrt{(l+1/2)^2-Z^2\alpha^2}=n-\sigma_l.\]
Where $\sigma_l=l+1/2-\sqrt{(l+1/2)^2-Z^2\alpha^2}$. And according
to (\ref{a5}), we have
\[\beta=\frac{2mZe_s{}^2}{\alpha'\hbar^2}=\frac{Ze_s{}^2}{\hbar}\left[\frac{m^2}{(m_0+m)|E'|}\right]^{1/2}\quad\mbox{or}\quad
(n-\sigma_l)^2=\frac{Z^2e_s{}^4}{\hbar^2}\frac{m^2}{(m_0+m)|E'|}.\]
Considering $m=m_0-|E'|/c^2$, then $|E'|$ satisfies the following
second-order algebraic equation:
\[\frac{Z^2\alpha^2+(n-\sigma_l)^2}{c^2}|E'|^2-2m_0[Z^2\alpha^2+(n-\sigma_l)^2]|E'|+Z^2\alpha^2m_0{}^2c^2=0.\]

Obviously, the expression of $|E'|$ obtained by solving the equation
is related to both $n=1,2,\ldots$ and $l=0,1,\ldots,n-1$, thus
$E=m_0c^2-|E'|$ can be denoted by $E_n$, then we have
\begin{eqnarray*}
  E_n &=& m_0c^2-|E'| \\
      &=&
      m_0c^2-m_0c^2\left(1\pm\frac{n-\sigma_l}{\sqrt{Z^2\alpha^2+(n-\sigma_l)^2}}\right)
      =\mp\frac{(n-\sigma_l)m_0c^2}{\sqrt{Z^2\alpha^2+(n-\sigma_l)^2}}.
\end{eqnarray*}

If simply taking the positive solution, then we have
\[E_n=\frac{m_0c^2}{\sqrt{1+Z^2\alpha^2/(n-\sigma_l)^2}}.\]

Therefore, for a hydrogen-like atom formed by pion capture in an
atomic nucleus, if we ignore the effects of nuclear motion, then its
relativistic energy level is
\begin{eqnarray}\label{a2}
  E_n &=& m_0c^2\left[1+\frac{Z^2\alpha^2}{\left(n-(l+1/2)+\sqrt{(l+1/2)^2-Z^2\alpha^2}\right)^2}\right]^{-1/2}\nonumber \\
      &=&
      m_0c^2\left[1-\frac{Z^2\alpha^2}{2n^2}-\frac{Z^4\alpha^4}{2n^4}\left(\frac{n}{l+1/2}-\frac{3}{4}\right)+\cdots\right].
\end{eqnarray}
Where $n=1,2,\ldots$ is the principal quantum number,
$l=0,1,\ldots,n-1$ is the angular quantum number, $m_0$ is the rest
mass of $\pi^-$, $Z$ is the atomic number, $\alpha$ is the fine
structure constant, $c$ is the speed of light in vacuum. The result
is completely the same as that of reference \cite{w}.

As $E_n=mc^2$, the system mass $m$ corresponding to the positive
solution is
\begin{equation}\label{a12}
    m=\frac{m_0}{\sqrt{1+Z^2\alpha^2/(n-\sigma_l)^2}}.
\end{equation}
Therefore, when the system is in a bound state, the system mass
takes discrete values.

According to (\ref{15}) and (\ref{16}), for a pionic hydrogen atom,
the solution of the Klein-Gordon equation has the probability
significance, further, we ought to obtain the relativistic wave
functions of such a system. Solving the radial equation (\ref{a3})
can come down to solving the equation (\ref{a7}), therefore, by
substituting the solutions $s=l+1-\sigma_l,\;\beta=n-\sigma_l$ of
the equation (\ref{a10}) into (\ref{a9}), we have
\[b_{\nu+1}=\frac{\nu+l+1-n}{(\nu+1-\sigma_l)(2l+2+\nu-\sigma_l)+Z^2\alpha^2}b_\nu.\]
\begin{eqnarray*}
  b_\nu &=& \frac{(l-n+\nu)(l-n+\nu-1)\cdots(l-n+2)(l-n+1)}{\prod^\nu_{k=1}[(k-\sigma_l)(2l+1+k-\sigma_l)+(Z\alpha)^2]}b_0 \\
    &=& \frac{(-1)^\nu(n-l-1)(n-l-2)\cdots(n-l-\nu)}{\prod^\nu_{k=1}[(k-\sigma_l)(2l+1+k-\sigma_l)]
    \prod^\nu_{k=1}\left(1+\frac{Z^2\alpha^2}{(k-\sigma_l)(2l+1+k-\sigma_l)}\right)}b_0 \\
    &=&
    \frac{(-1)^\nu(n-l-1)!}{(n-l-1-\nu)!\prod^\nu_{k=1}[(k-\sigma_l)(2l+1+k-\sigma_l)]\eta(l,\nu)}b_0.
\end{eqnarray*}
Where $b_0$ is a constant:
\[b_0=-\frac{[(n+l)!\,]^2}{(n-l-1)!\,\Gamma(1-\sigma_l)\Gamma(2l+2-\sigma_l)}\frac{N_{nl}}{\alpha'},\]
where $\alpha'$ is defined by (\ref{a5}), thus we have
\[b_\nu=\frac{(-1)^{\nu+1}[(n+l)!\,]^2}{(n-l-1-\nu)!\,\Gamma(\nu+1-\sigma_l)\Gamma(2l+2+\nu-\sigma_l)\eta(l,\nu)}\frac{N_{nl}}{\alpha'}.\]
Substituting it into (\ref{a8}), then
$R=u(\rho)/r=e^{-\rho/2}f(\rho)/r$ is definitely related to both $n$
and $l$, denoted by $R_{nl}(r)$, thus we have
\[R_{nl}(r)=N_{nl}e^{-\rho/2}\sum^{n-l-1}_{\nu=0}
\frac{(-1)^{\nu+1}[(n+l)!\,]^2\rho^{l-\sigma_l+\nu}}{(n-l-1-\nu)!\,\Gamma(\nu+1-\sigma_l)\,\Gamma(2l+2+\nu-\sigma_l)\,\eta(l,\nu)}.\]
According to (\ref{a5}), we have
\[\rho=\alpha'r=\frac{2mZe_s{}^2}{\hbar^2\beta}r=\frac{2mZe_s{}^2}{(n-\sigma_l)\hbar^2}r=\frac{2Z}{(n-\sigma_l)a_0}r,
\quad{a_0}=\frac{\hbar^2}{me_s{}^2}.\]

Therefore, for a hydrogen-like atom formed by pion capture in an
atomic nucleus, if we ignore the effects of nuclear motion, then its
relativistic stationary state wave function is
\begin{equation}\label{a13}
    \Psi(\mathbf{r},t)=\psi_{nlm}(r,\theta,\varphi)e^{-iE_nt/\hbar}=R_{nl}(r)Y_{lm}(\theta,\varphi)e^{-iE_nt/\hbar}.
\end{equation}
Where $Y_{lm}(\theta,\varphi)$ is the normalized spherical
harmonics. Relativistic effects do not cause the change of the
angular wave function $Y_{lm}(\theta,\varphi)$, but the change of
the radial wave function $R_{nl}(r)$, i.e.
\begin{equation}\label{a14}
R_{nl}(r) =
N_{nl}e^{-\frac{Z}{(n-\sigma_l)a_0}r}\left(\frac{2Z}{(n-\sigma_l)a_0}r\right)^{l-\sigma_l}
L^{2l+1-\sigma_l}_{n+l}\left(\frac{2Z}{(n-\sigma_l)a_0}r\right).
\end{equation}
\[\sigma_l=l+\frac{1}{2}-\left[\left(l+\frac{1}{2}\right)^2-Z^2\alpha^2\right]^{1/2}=
\sum^\infty_{k=1}\frac{2^{k-1}(2k-3)!!}{k!(2l+1)^{2k-1}}(Z\alpha)^{2k}.\]
Where $n$ is the principal quantum number, $l$ is the angular
quantum number, $m$ is the magnetic quantum number (do not confuse
it with the system mass $m$), $Z$ is the atomic number,
$\alpha=e_s{}^2/\hbar{c}$ is the fine structure constant, and
\[a_0=\frac{\hbar^2}{me_s{}^2},\qquad{m}=m_0\left[1+\frac{Z^2\alpha^2}{(n-\sigma_l)^2}\right]^{-1/2}.\]
Where $m$ is the system mass, and $m_0$ is the rest mass of a
particle. The expression of $L^{2l+1-\sigma_l}_{n+l}(\rho)$ is
\begin{equation}\label{a15}
L^{2l+1-\sigma_l}_{n+l}(\rho) = \sum^{n-l-1}_{\nu=0}
\frac{(-1)^{\nu+1}[(n+l)!\,]^2\rho^\nu}{(n-l-1-\nu)!\,\Gamma(2l+\nu+2-\sigma_l)\Gamma(\nu+1-\sigma_l)\,\eta(l,\nu)},
\end{equation}
which is called the relativistic associated Laguerre polynomials,
where
\[\eta(l,\nu)=\prod^\nu_{k=1}\left(1+\frac{Z^2\alpha^2}{(k-\sigma_l)(2l+1+k-\sigma_l)}\right).\]
$R_{nl}(r)$ is also normalized, and $N_{nl}$ is the normalized
constant.

According to (\ref{15}) and (\ref{16}), by using the Gauss's theorem
in vector analysis, the relativistic stationary state wave function
$\Psi(\mathbf{r},t)$ satisfies the normalization condition:
$\int_\infty{w}(\mathbf{r},t)\,d\tau=\int_\infty\Psi^*(\mathbf{r},t)\Psi(\mathbf{r},t)\,d\tau=1$,
which can also be expressed as
\[\int^\infty_0R^2_{nl}(r)r^2\,dr=1\quad\mbox{and}\quad
\int^\pi_0\int^{2\pi}_0Y^*_{lm}(\theta,\varphi)Y_{lm}(\theta,\varphi)\sin\theta\,d\theta\,d\varphi=1.\]
Thus the normalized constant $N_{nl}$ can be determined.

The value of the fine structure constant is very small
($\alpha\approx1/137$), for $Z=1$, in the non-relativistic limit,
$\eta(l,\nu)\rightarrow1,\;\sigma_l\rightarrow0$, Thus
$L^{2l+1-\sigma_l}_{n+l}(\rho)\rightarrow{L^{2l+1}_{n+l}(\rho)}$.
$L^{2l+1}_{n+l}(\rho)$ is the associated Laguerre polynomials.

\section{The Klein-Gordon Equation}

Now, let us prove that the wave equations (\ref{7}) and (\ref{9})
are relativistic. Assuming there are two equations, one of which is
relativistic, while it is unknown whether or not the other is
relativistic, if we can prove these two equations have the same
solution, then the other is also relativistic. According to this
thinking, if we directly quantize (\ref{5}), then
$\psi(\mathbf{r})$, which is the main part of the stationary state
wave function, satisfies the following wave equation:
\[\left(E'+m_0c^2-U(\mathbf{r})\right)^2\psi=m_0{}^2c^4\psi-c^2\hbar^2\nabla^2\psi,\]
it has the same solution as (\ref{10}). We can say they are
equivalent, and it can be equivalently expressed as:
\begin{equation}\label{10''}
\left(E-U(\mathbf{r})\right)^2\psi=m_0{}^2c^4\psi-c^2\hbar^2\nabla^2\psi.
\end{equation}
Further, we have the more general form of this equation
\begin{equation}\label{10'}
\left(i\hbar\frac{\partial}{\partial{t}}-U(\mathbf{r})\right)^2\Psi=m_0{}^2c^4\Psi-c^2\hbar^2\nabla^2\Psi.
\end{equation}

Under strong coupling conditions, the relativistic effects of moving
particles in a potential field become significant, therefore,
seeking for the exact solution of the Klein-Gordon equation or the
Dirac equation in a typical potential field has drawn increasing
attention in the last few years. Let $m_0$ be the rest mass, and $E$
be the total energy containing the rest energy, according to
reference \cite{daf}-\cite{cls}, the Klein-Gordon equation with the
vector potential $U(\mathbf{r})$ and the scalar potential
$S(\mathbf{r})$ in spherical coordinates is:
\begin{equation}\label{y1}
[-\hbar^2c^2\nabla^2+(m_0c^2+S(\mathbf{r}))^2]\psi(\mathbf{r})=[E-U(\mathbf{r})]^2\psi(\mathbf{r}).
\end{equation}
When $S(\mathbf{r})=0$, we have (\ref{10''}). In other words, by
introducing the concept of system mass, the Klein-Gordon equation
can be expressed in a better form when compared with the
non-relativistic limit, making it easier to transplant the concepts
and methods of the non-relativistic quantum mechanics into the
relativistic quantum mechanics.

If we introduce the concept of system mass, (\ref{y1}) can be
expressed as

\parbox{11cm}{\begin{eqnarray*}\label{y2}
E'\psi&=&-\,\frac{\hbar^2}{m_0+m}\nabla^2\psi+\frac{2m}{m_0+m}U\psi-\frac{U^2}{(m_0+m)c^2}\psi\\
&&+\,\frac{2m_0}{m_0+m}S\psi+\frac{S^2}{(m_0+m)c^2}\psi.
\end{eqnarray*}}\hfill\parbox{1cm}{\begin{eqnarray}\end{eqnarray}}
Where $m$ is the system mass, and $E'=E-m_0c^2$, $E=mc^2$.
(\ref{y2}) formally remains similar with the Schr\"{o}dinger
equation, and the solution of the equation also has the probability
significance. Now, let us seek for the relation between the
solutions of (\ref{10'}) in two different inertial frames of
reference. Let $K$ and $K'$ be two inertial frames of reference, and
their coordinate axes be correspondingly parallel, $K'$ be moving at
a speed of $v$ along the positive x-axis relative to $K$, and the
coincident moment of $O$ and $O'$ be the starting time point. In the
inertial frame of reference $K$, particles are moving in the
potential field $U$, while in the inertial frame of reference $K'$,
particles are moving in the potential field $U'$. We use the
following transformation relation as the characteristic equation for
a single-particle system, i.e. the Lorentz transformation from $K$
to $K'$ is

\parbox{10cm}{\begin{eqnarray*}\label{13}
p_x{}'&=&\frac{1}{\sqrt{1-\beta^2}}p_x-\frac{v/c^2}{\sqrt{1-\beta^2}}(E-U)\\
p_y{}'&=&p_y\\p_z{}'&=&p_z\\
E'-U'&=&\frac{1}{\sqrt{1-\beta^2}}(E-U)-\frac{v}{\sqrt{1-\beta^2}}p_x.
\end{eqnarray*}}\hfill\parbox{1cm}{\begin{eqnarray}\end{eqnarray}}

The Lorentz transformation from $K'$ to $K$ is

\parbox{11cm}{\begin{eqnarray*}\label{14}
p_x&=&\frac{1}{\sqrt{1-\beta^2}}p_x{}'+\frac{v/c^2}{\sqrt{1-\beta^2}}(E'-U')\\
p_y&=&p_y{}'\\p_z&=&p_z{}'\\
E-U&=&\frac{1}{\sqrt{1-\beta^2}}(E'-U')+\frac{v}{\sqrt{1-\beta^2}}p_x{}'.
\end{eqnarray*}}\hfill\parbox{1cm}{\begin{eqnarray}\end{eqnarray}}
Here $E'$ is the total energy of the particle moving in the
potential field $U'$, which is different from the $E'$ in Definition
2.

In the base function
$\Psi_p(\mathbf{p},t)=A\exp[-i(Et-\mathbf{p\cdot{r}})/\hbar]$,
$Et-\mathbf{p\cdot{r}}$ is the phase of a plane wave, when
transformed from $K$ to $K'$, $E't'-\mathbf{p'\cdot{r}'}$ and
$Et-\mathbf{p\cdot{r}}$ differ by at most an integral multiple of
$2\pi$. By the periodicity of phase, we have $Et-\mathbf{p\cdot{r}}$
is a relativistic invariant. Therefore, by choosing $A'$ properly,
we have
\[A\exp[-i(Et-\mathbf{p\cdot{r}})/\hbar]=A'\exp[-i(E'{t}'-\mathbf{p'\cdot{r}'})/\hbar]\quad
\mbox{or}\quad\Psi_p(\mathbf{r},t)=\Psi_p{}'(\mathbf{r'},t').\]

Thus if we multiply the left-hand side of the characteristic
equation (\ref{13}) by $\Psi_p{}'(\mathbf{r'},t')$, and multiply the
right-hand side of the equation by $\Psi_p(\mathbf{r},t)$, the
equation is still tenable. By using (\ref{2}) and similar
differential relations, we have the following set of relations:
\begin{eqnarray*}
-i\hbar\frac{\partial\Psi_p{}'}{\partial{x'}}&=&-\frac{i\hbar}{\sqrt{1-\beta^2}}\left(\frac{\partial\Psi_p}{\partial{x}}+\frac{v}{c^2}\frac{\partial\Psi_p}{\partial{t}}\right)+\frac{v/c^2}{\sqrt{1-\beta^2}}U\Psi_p\\
-i\hbar\frac{\partial\Psi_p{}'}{\partial{y'}}&=&-i\hbar\frac{\partial\Psi_p}{\partial{y}}\\
-i\hbar\frac{\partial\Psi_p{}'}{\partial{z'}}&=&-i\hbar\frac{\partial\Psi_p}{\partial{z}}\\
-U'\Psi_p{}'+i\hbar\frac{\partial\Psi_p{}'}{\partial{t'}}&=&-\frac{1}{\sqrt{1-\beta^2}}U\Psi_p+\frac{i\hbar}{\sqrt{1-\beta^2}}\left(\frac{\partial\Psi_p}{\partial{t}}+v\frac{\partial\Psi_p}{\partial{x}}\right).
\end{eqnarray*}
By using (\ref{1a}), we have:

Let $K$ and $K'$ be two inertial frames of reference, and their
coordinate axes be correspondingly parallel, $K'$ be moving at a
speed of $v$ along the positive x-axis relative to $K$, and the
coincident moment of $O$ and $O'$ be the starting time point, then
for a single-particle system, the Lorentz transformation of wave
functions from $K$ to $K'$ is

\parbox{10cm}{\begin{eqnarray*}\label{11}
-i\hbar\frac{\partial\Psi'}{\partial{x'}}&=&-\frac{i\hbar}{\sqrt{1-\beta^2}}\left(\frac{\partial\Psi}{\partial{x}}+\frac{v}{c^2}\frac{\partial\Psi}{\partial{t}}\right)+\frac{v/c^2}{\sqrt{1-\beta^2}}U\Psi\\
-i\hbar\frac{\partial\Psi'}{\partial{y'}}&=&-i\hbar\frac{\partial\Psi}{\partial{y}}\\
-i\hbar\frac{\partial\Psi'}{\partial{z'}}&=&-i\hbar\frac{\partial\Psi}{\partial{z}}\\
-U'\Psi'+i\hbar\frac{\partial\Psi'}{\partial{t'}}&=&-\frac{1}{\sqrt{1-\beta^2}}U\Psi+\frac{i\hbar}{\sqrt{1-\beta^2}}\left(\frac{\partial\Psi}{\partial{t}}+v\frac{\partial\Psi}{\partial{x}}\right).
\end{eqnarray*}}\hfill\parbox{1cm}{\begin{eqnarray}\end{eqnarray}}

Similarly, the Lorentz transformation of wave functions from $K'$ to
$K$ is

\parbox{11cm}{\begin{eqnarray*}\label{12}
-i\hbar\frac{\partial\Psi}{\partial{x}}&=&-\frac{i\hbar}{\sqrt{1-\beta^2}}\left(\frac{\partial\Psi'}{\partial{x'}}-\frac{v}{c^2}\frac{\partial\Psi'}{\partial{t'}}\right)-\frac{v/c^2}{\sqrt{1-\beta^2}}U'\Psi'\\
-i\hbar\frac{\partial\Psi}{\partial{y}}&=&-i\hbar\frac{\partial\Psi'}{\partial{y'}}\\
-i\hbar\frac{\partial\Psi}{\partial{z}}&=&-i\hbar\frac{\partial\Psi'}{\partial{z'}}\\
-U\Psi+i\hbar\frac{\partial\Psi}{\partial{t}}&=&-\frac{1}{\sqrt{1-\beta^2}}U'\Psi'+\frac{i\hbar}{\sqrt{1-\beta^2}}\left(\frac{\partial\Psi'}{\partial{t'}}-v\frac{\partial\Psi'}{\partial{x'}}\right).
\end{eqnarray*}}
\hfill\parbox{1cm}{\begin{eqnarray}\end{eqnarray}} Where
$\beta=v/c$, $U$ is the potential filed in $K$, and $U'$ is the
potential filed in $K'$.

\section{Conclusion}

In summary, by introducing Definition 1-5, using (\ref{1a}) and
assuming that the relativistic kinetic expression is tenable in a
larger sense, the Klein-Gordon equation in a potential field
(\ref{10'}) is derived more strictly, and we also prove that its
stationary solutions ($\psi(\mathbf{r})e^{-iEt/\hbar}$) satisfy the
equation (\ref{7}), thus they have the probability significance.

Under the condition that the scalar and vector potentials are equal,
$S(\mathbf{r})=U(\mathbf{r})$, (\ref{y2}) degenerates into
\begin{equation}\label{y3}
E'\psi=-\,\frac{\hbar^2}{m_0+m}\nabla^2\psi+2U\psi,
\end{equation}
which is formally similar to the Schr\"{o}dinger equation.
Therefore, the following conjecture brought out in reference
\cite{cls} is correct, i.e. if the bound states of the
Schr\"{o}dinger equation in a certain potential field can be solved
exactly, then under the condition that the scalar and vector
potentials are equal, the bound states of the Klein-Gordon equation
in the same potential field can also be solved exactly by using the
same method. Reference \cite{zcd}-\cite{cls} have given several
exact bound state solutions of the Klein-Gordon equation and Dirac
equation in a non-central potential field with equal scalar and
vector potentials. And under the condition that the scalar and
vector Manning-Rosen potentials are equal, the bound state solutions
of the s-wave Klein-Gordon equation and Dirac equation are
respectively derived in reference \cite{zw}. By using (\ref{y3}),
these results can be derived much easier.

According to reference \cite{cls2}, the radial Schr\"{o}dinger
equation with the Hulth\'{e}n potential is
\begin{equation}\label{y4}
\frac{d^2u(r)}{dr^2}+\frac{2m_0}{\hbar^2}\left(E'+Ze^2\lambda\frac{e^{-\lambda{r}}}{1-e^{-\lambda{r}}}-\frac{l(l+1)\hbar^2}{2m_0r^2}\right)u(r)=0.
\end{equation}
According to (\ref{y3}), under the condition that the scalar and
vector potentials are equal, the radial Klein-Gordon equation with
the Hulth\'{e}n potential is
\begin{equation}\label{y5}
\frac{d^2u(r)}{dr^2}+\frac{m_0+m}{\hbar^2}\left(E'+2Ze^2\lambda\frac{e^{-\lambda{r}}}{1-e^{-\lambda{r}}}-\frac{l(l+1)\hbar^2}{(m_0+m)r^2}\right)u(r)=0.
\end{equation}
Obviously, solving this equation does not increase any difficulty,
thus the results of the scattering states of the Hulth\'{e}n
potential in reference \cite{cls2} can be extended from the
Schr\"{o}dinger equation to the Klein-Gordon equation with equal
scalar and vector potentials. Therefore, we have also come to the
conclusion that, if the scattering states of the Schr\"{o}dinger
equation in a certain potential field can be solved exactly, then
under the condition that the scalar and vector potentials are equal,
the scattering states of the Klein-Gordon equation in the same
potential field can also be solved exactly by using the same method.

\end{document}